\title{ ~~\\ Around Pelik\'an's conjecture on very odd sequences}
\author{Pieter Moree and Patrick Sol\'e}
\documentclass[12pt]{article}
\usepackage{amssymb, latexsym, amsfonts}
\def\@ptsize{2}
\newtheorem{Thm}{Theorem}
\newtheorem{Conj}{Conjecture}
\newtheorem{Lem}{Lemma}

\newtheorem{Cor}{Corollary}

\newtheorem{Prop}{Proposition}
\newcommand{\qed}{\hfill $\Box$}

\begin{document}
\date{}
\maketitle
{\def\thefootnote{}
}
\begin{abstract}
\noindent Very odd sequences were introduced in 1973 by Pelik\'an who
conjectured that
there were none of length $\ge 5.$ This conjecture was disproved first by
 MacWilliams and Odlyzko \cite{MO} in 1977 and then  by two different sets of authors
in 1992 \cite{A}, 1995 \cite{IW}. 
We give connections with duadic codes, cyclic difference sets, levels (Stufen) of
cyclotomic fields,
and derive some new asymptotic results on the length of very odd sequences and
the number of such sequences of a given length.
\end{abstract}
\section{Introduction}
For a given natural number $n$ fix integers $a_i$ with
$a_i\in \{0,1\}$ for $1\le i\le n$. 
Put $A_k=\sum_{i=1}^{n-k}a_ia_{i+k}$ for $0\le k\le n-1$.
We say that $a_1\ldots a_n$ is a {\it very odd sequence} if $A_k$ is odd for
$1\le k\le n$.
By $S(n)$ we denote the number of very odd sequences of length $n$.
Pelik\'an \cite{P} conjectured that very odd sequences
of length $n\ge 5$ do not exist. However,
the sequence $101011100011$ (for which the corresponding $A_k$'s are
$7,3,3,1,3,3,3,1,1,1,1,1$) disproves this conjecture. Alles \cite{A}
showed that if $S(n)>0$, then also $S(7n-3)>0$ and thus showed that there
are infinitely many counterexamples to Pelik\'an's conjecture. In his
note \cite{A} Alles raised two questions:\\
(1) Does $S(n)>0$ imply $n\equiv 0,1({\rm mod~}4)$?\\
(2) Does $S(n)>0$ imply $S(n)=2^k$ for some integer $k$?\\ 
These two questions, which also appear as 
unsolved problem E38 in \cite{G}, were positively answered in \cite{IW}.\\
\indent All of the above was, however, already known much earlier \cite{MO}. 
For $a$ and $b$ coprime
integers, let ord$_a(b)$ denote the smallest positive integer $r$ such
that $a^r\equiv 1({\rm mod~}b)$. Let ${\cal P}=\{7,23,31,47,71,73,79,\ldots\}$ denote the set of odd primes
$p$ for which ord$_2(p)$ is odd (throughout this paper the letter $p$ will be used
to denote primes).
MacWilliams
and Odlyzko proved that $S(n)>0$ iff $2n-1$ is composed only 
from primes in $\cal P$. (Alternatively we can formulate this as
$S(n)>0$ iff the order of 2 modulo $2n-1$ is odd.) Since $7\in \cal P$,
the result of Alles that $S(n)>0$ implies $S(7n-3)>0$, immediately follows.
Using the supplementary law of quadratic reciprocity it is easily
seen, as already noticed in \cite{MO}, that $p\in \cal P$ implies $p\equiv \pm 1({\rm mod~}8)$.
Thus if $S(n)>0$, then $2n-1\equiv \pm 1({\rm mod~}8)$ and from this it follows that
the answer to question (1) is affirmative. In \cite{MO} question (2) was answered
in the affirmative also; to be more precise it was shown that if $S(n)>0$ then $X^{2n-1}+1$
decomposes over $\mathbb F_2$ into an odd number of distinct irreducible factors, say
$2h+1$ irreducible factors in total and that, moreover, $S(n)=2^h$. R. van der Veen 
and E. Nijhuis \cite{VN} described and implemented an algorithm to
determine all very odd sequences of a given length. Moreover, they showed that, for $n>1$,
a very odd sequence is never periodic.\\ 
\indent Let $N(x)$ denote the number of $n\le x$ for which very odd sequences of length $n$ exist.
Let Li$(x)$ denote the logarithmic integral. By \cite[Theorem 2]{M} 
it follows that
\begin{equation}
\label{scherper1}
{\cal P}(x)
={7\over 24}{\rm Li}(x)+O\left({x(\log \log x)^4\over \log^3 x}\right),
\end{equation}
where ${\cal P}(x)$ denotes the number of primes $p\le x$ in $\cal P$.
Using Theorem 4 of \cite{M} it then follows, that there exists 
a $c_0>0$ such that
\begin{equation}
\label{scherper2}
N(x)=c_0{x\over \log^{17/24}x}\left(1+O\left({(\log \log x)^5\over 
\log x}\right)\right).
\end{equation}
Since $\lim_{x\rightarrow \infty}N(x)/x=0$,
it follows in particular that Pelik\'an's conjecture holds
true for almost all integers $n$. The estimates (\ref{scherper1}) 
and (\ref{scherper2}) sharpen the assertions ${\cal P}(x)\sim {7\over 24}
{x\over \log x}$, 
respectively $N(x)=o(x)$ made in \cite{MO}.\\
\indent The analysis of very odd sequences and that
of Ducci sequences \cite{BLM} shows a certain analogy (in the latter case the factorization of
$(1+X)^n+1$ over $\mathbb F_2[X]$ plays an important role and there is also  a link
with Artin's primitive root conjecture).\\
\indent In Sections 2 and 3 the connection between very odd sequences and coding theory, 
respectively the Stufe (level) of cyclotomic fields is considered. 
These sections have a partly survey nature and can be read independently from the remaining sections.
In Section 4 a formula
for $S(n)$ is derived (Proposition \ref{Sformula}). In Section 5 the value distribution
of $S(n)$ is then considered, using some results related to Artin's primitive root
conjecture (Theorem \ref{stellinkje} and Theorem \ref{stellinkje2}). It is shown that
under the Generalized Riemann Hypothesis (GRH) the preimage of the value $2^e$ is 
infinite, provided we can find just one integer $n$ of a certain type satisfying
$S(n)=2^e$ (Theorem \ref{eenvoorallen}). This leads us then into a query of finding
such integers and motivates the introduction of the notion of solution tableaux (Section 
6.1). In the final section it is shown that certain values of $S(n)$ are assumed much more 
infrequently than others.

\section{Codes and Pelik\'an's conjecture}
For coding theoretic terminology, we
    refer to \cite{MS77}.
\begin{Lem}
We have $S(n)>0$ iff there exists a $[2n,n]$ extended binary cyclic code.
In particular, there are very odd sequences of length ${p+1\over 2}$
for all odd  primes $p \equiv \pm 1({\rm mod~}8),$ and of length $2^m$
 for all even integers $m\ge 2.$
\end{Lem}
{\it Proof}. Write $b_k:=a_{k+1}$ and $b(X):=\sum_{k=0}^{n-1} b_k X^k$,
 with  $X$ an indeterminate and consider $b(X)$ as
 an element of $\mathbb F_2[X]$. 
 We say $b(X)$ is the polynomial associated to $a_1\ldots a_n$.
 It was already observed in \cite{MO} that 
$a_1\ldots a_n$ is a very odd sequence iff the polynomial identity
$X^{2n-1}+1=(X+1)b(X)b^*(X),$
the $*$ denoting reciprocation, holds true.
Using this observation we see that the cyclic code of length $2n-1$ and generator
$b(X)$ is self-orthogonal of dimension $n.$ By adding an overall
parity-check we obtain a self-dual code of length $2n$.\\
\indent The two infinite families of lengths
correspond to, respectively, quadratic residue codes attached 
to the prime $p$ and Reed-Muller codes RM$({m\over 2},m+1)$ (see 
\cite{MS77}).\qed\\

\noindent We denote by $C(n)$ the number of nonequivalent codes of length
$2n$ obtained by the preceding lemma. 
The celebrated $[24,12,8]$ Golay code
arises on taking $n=12$. More generally, we obtain duadic codes with 
multiplier $-1$ which have received much attention in the last 
twenty years \cite{p1,p2}. In the special case of $n$ a multiple 
of $8,$ these codes are, furthermore, doubly even. There are 
generalizations over $\mathbb F_4$ \cite{p3}, $\mathbb F_q$ \cite{Smi,Smi2} and
over rings \cite{LaS,LiS}.
The following connection with difference sets is proved
differently in \cite[Thm 6.2.1.]{Smi2} and in a very general setting in \cite{Ru}.
\begin{Lem}
If there is a $2-(N,K,\lambda)$ cyclic difference set with $K$ and $\lambda$
both 
odd, then there is a very odd sequence of length $(N+1)/2.$
\end{Lem}
{\it Proof}.
We use the ring morphism
(given by reduction mod 2) between the algebra $\mathbb F_2[X]/(X^N-1),$
where binary cyclic codes
of length $N$ live and the group algebra
$\mathbb Z[C_N]$ which occurs in the following
characterization \cite[Lemma 3.2, p. 312]{BLJ} of difference sets of the cyclic
group $C_N:$ $$DD^{-1}=K-\lambda+\lambda C_N.$$
The difference set $D$ is then the set of 
exponents of $X$ occurring in the polynomial $b(X)$ of the preceding proof. \qed\\

\noindent The following construction generalizes Alles's \cite{A}, who took
$m=4$.
Let $a$ (resp. $b$) denote a very odd sequence of length $n$
(resp. $m$).
Let $a(X)$ (resp. $b(X)$) denote their associated polynomials in
$\mathbb F_2[X]$.
By $a\otimes b$ we shall mean the sequence corresponding to the polynomial
$a(X)b(X^{2n-1})$. Explicitly this amounts to the sequence obtained by taking $b$, replacing each
$b_i\in b$ by a block of $n$ zeros (if $b_i=0$) or by the
sequence $a$ (if $b_i=1$), and
finally inserting a sequence of $n-1$ zeros between each of these blocks.
The following result is given without a proof in \cite{IW}.
\begin{Lem}
If $a$ and $b$ are very odd sequences of respective lengths $n$ and $m,$
then $a\otimes b$ is a very odd sequence of length $2mn-n-m+1.$
If ${a}'$ and ${b'}$ are very odd sequences of lengths $n$ and $m$, then
$a\otimes b=a'\otimes b'$ implies $a=a'$ and $b=b'$.
\end{Lem}
{\it Proof}. We know that, by hypothesis, $Y^{2m-1}+1=(Y+1)b(Y)b^*(Y).$
Letting $Y=X^{2n-1}$ and using the hypothesis on $a$, that is
$X^{2n-1}+1=(X+1)a(X)a^*(X),$
we see that $c(X):=a(X)b(X^{2n-1}),$ satisfies
$X^{2N-1}+1=(X+1)c(X)c^*(X),$
for $2N-1=(2n-1)(2m-1),$ i.e. $N=2mn-n-m+1.$\\
\indent Conversely, writing c(X)=A(X)B(X), with  $A,B$ such that
\begin{itemize}
\item the order of the roots of $A$ divide $2n-1$
\item the order of the roots of $B$ does {\it not}  divide $2n-1$
\end{itemize}
we see that such a factoring is unique. We claim that $A=a$ and
$B(X)=b(X^{2n-1})$ is of that shape. Indeed, if $\beta^{2n-1}=1$ and
$b(\beta^{2n-1})=0$ we have that $1$ is a root of $b(X)$, contradicting the
definition of $b.$ Alternatively note that the second part of the assertion
follows from the explicit construction of $a\otimes b$.
\qed
\begin{Cor}
\label{wiseguys}
We have $S(2mn-n-m+1)\ge S(m)S(n)$.   
\end{Cor}

\section{Connection with the Stufe}
Let $K$ be a field. Then the Stufe (or level) of $K$, $s(K)$, is defined as the smallest number $s$ 
(if this
exists) such that $-1=\alpha_1^2+\cdots+\alpha_s^2$ with all $\alpha_i\in K$. 
If the Stufe does not exist, it is not difficult to see that there exists an order $\le$ on $K$ that
is compatible with the field operations (i.e. $K$ is orderable).
Pfister \cite{pfis}
proved that the Stufe of any field, if it exists, is a power of two. For $m\ge 1$, let
$K_m=\mathbb Q(e^{2\pi i/m})$. Hilbert proved that $s(K_m)\le 4$ for $m\ge 3$. Moser 
\cite{Mos1, Mos2} and, independently, Fein
et al. \cite{FGS} proved that $s(K_{2m-1})=2$ iff ord$_2(2m-1)$ is even.
We thus obtain the following lemma.
\begin{Lem}
\label{Stufe}
We have $S(n)>0$ iff $s(K_{2n-1})=4$.
\end{Lem}
In the same vein Ji \cite{Ji} proved in case $m>1$ is odd that every {\it algebraic integer} in $K_m$ can
be expressed as a sum of three integral squares iff $ord_2(m)$ is even.\\
\indent Note that if $4|m$, then clearly $s(K_m)=1$ and for $m$ odd we have $K_{2m}=K_m$, thus
once we know $s(K_m)$ for every odd $m$, we know the Stufe of every cyclotomic field. 
Lemma \ref{Stufe} then shows
that knowing for which $n$ we 
have $S(n)>0$ is equivalent with knowing the Stufe 
for every cyclotomic field. Interestingly enough, knowing the Stufe of every imaginary quadratic
number field is equivalent with knowing which integers can be 
represented as a sum
of three integer squares (Gauss' famous three-squares theorem), see 
\cite{Sm}.  
Let $St_4(x)$ count the number of $m\le x$ such that $s(K_{2m-1})=4$.
Then from Lemma \ref{Stufe} and (\ref{scherper2}) we infer that 
$$St_4(x)=N(x)=c_0{x\over \log^{17/24}x}\left(1+O\left({(\log \log x)^5\over \log x}\right)\right).$$

\section{An explicit formula for $S(n)$}
In the introduction we already remarked that if $S(n)>0$, then $X^{2n-1}+1$ decomposes over 
$\mathbb F_2$ into an odd number of disctinct irreducible factors, say $2h+1$ factors in total and
that, moreover, $S(n)=2^h$. 
Using this it is not difficult to derive an explicit formula for $S(n)$ 
(Proposition \ref{Sformula}). To this end we first study
the factorization of $X^n-1$ into irreducibles 
over $\mathbb F_q$, with $q$ a power of a prime
$p$. (Note that $X^n-1$ and $X^n+1$ represent the same polynomial in $\mathbb F_2[X]$.)
If $n=mp^e$, with $p\nmid n$,
then since $X^n-1=(X^m-1)^{p^e}$ over $\mathbb F_q$, we can reduce to the case where
$(n,q)=1$. Then we have the following result.
\begin{Lem}
\label{intoirreducibles}
Let $q$ be the order of a finite field. If $(n,q)=1$, then $X^n-1$ factors
into 
\begin{equation}
\label{iformula}
i_q(n)=\sum_{d|n}{\varphi(d)\over {\rm ord}_q(d)}
\end{equation}
distinct irreducible factors over $\mathbb F_q$. 
\end{Lem}
{\it Proof}. The assumption $(n,q)=1$ ensures that the irreducibles will be
distinct. 
Let $\Phi_d(X)$ denote the cyclotomic polynomial of degree $d$. We have
$X^n-1=\prod_{d|n}\Phi_d(X)$, cf. \cite[Theorem 2.45]{LN}. It is
not difficult to show \cite[Theorem 2.47]{LN} that the polynomial $\Phi_d(X)$
decomposes into $\phi(d)/{\rm ord}_q(d)$ irreducibles over $\mathbb F_q$ each having
degree ${\rm ord}_q(d)$. On
combining these results, the proof is then completed. \qed\\

\noindent The above formula for $i_q(n)$ also arises in some other mathematical contexts,
see e.g. \cite{C,R,U,VS}. In particular we like to recall the 
following important result due to Ulmer \cite{U}.
\begin{Thm}
Let $p$ be a prime, $n$ a positive integer, and $d$ a divisor of $p^n+1$ that
is coprime with $6$. Let $q$ be a power of $p$ and let $E$ be the elliptic curve
over $\mathbb F_q(t)$ defined by $y^2+xy=x^3-t^d$. Then the $j$-invariant of $E$ is
not in $\mathbb F_q$, the conjecture of Birch and Swinnerton-Dyer holds for $E$, and
the rank of $E(\mathbb F_q(t))$ equals $i_q(d)$.
\end{Thm}
We suspect that some of the techniques of this paper can be used to study the
value distribution of the ranks in the latter result (for fixed $q$); a question
that seems of some importance.\\

\noindent The function $i_q$ is neither multiplicative nor additive, but we
can still prove something in this direction (Proposition \ref{inequality}). The proof makes use of
the following lemma.  
\begin{Lem}
\label{deelbaarheid} 
For $q\nmid n$, let $r_q(n)={\varphi(n)/{\rm ord}_q(n)}$.\\
\noindent {\rm 1)} Suppose that
$\delta|d$ and $(d,q)=1$, then $r_q(\delta)|r_q(d)$.\\
{\rm 2)} If $r_q(p)=r_q(p^2)$, then $r_q(p^e)=r_q(p)$ for every $e\ge 1$.
\end{Lem}
{\it Proof}. 1) The natural projection of multiplicative groups
$(\mathbb Z/d\mathbb Z)^*\rightarrow  (\mathbb Z/\delta\mathbb Z)^*$ gives
rise to the projection $(\mathbb Z/d\mathbb Z)^*/\langle q\rangle 
\rightarrow (\mathbb Z/\delta\mathbb Z)^*/\langle q\rangle$, and
so $r_q(\delta)$ divides $r_q(d)$ as claimed.\\
2) It is a well-known, and easy to prove, result in elementary number
theory that if ord$_q(p)=m$ and ord$_q(p^2)=pm$, then 
${\rm ord}_q(p^e)=p^{e-1}m$ for $e\ge 1$. \qed\\

\noindent {\tt Remark}. Note
that $r_q(n)=[(\mathbb Z/n\mathbb Z)^*:\langle q\rangle]$. This quantity
is sometimes called the (residual) index of $q$ in  $(\mathbb Z/n\mathbb Z)^*$.\\  

\noindent {\tt Remark}. Let $f(X)\in \mathbb Z[X]$ be an 
irreducible monic polynomial
over $\mathbb Q$. A celebrated result of Dedekind (see 
e.g. \cite[Theorem 4.12]{N}) relates the factorization of
$f(X)$ over $\mathbb F_p$ to the factorization of the ideal $(p)$ into prime ideals
in the ring of integers $\cal O$ of the quotient field $\mathbb Q[X]/f(X)$ (to 
each irreducible $f_i(X)$  of $f(X)$ over $\mathbb F_p$ corresponds a prime ideal 
lying over $(p)$ in $\cal O$ of degree deg$(f_i)$). When we apply
this with $f(X)=\Phi_d(X)$ and $p\nmid d$, we see that
in $\mathbb Z[\zeta_d]$, the ring of integers of $\mathbb Q(\zeta_d)$, the
prime ideal $(p)$ factorizes as $(p)={\cal P}_1\cdots {\cal P}_g$, where
each ${\cal P}_i$ has degree ord$_p(d)$ and $g=r_p(d)$ (see e.g. \cite[Theorem 4.16]{N}).
This interpretation of $r_p(d)$ together with some basic facts from algebraic
number theory gives another proof of Lemma \ref{deelbaarheid} in 
case $q=p$. {}From the above remarks it follows 
that $i_p(n)$, for $p\nmid n$, counts the total number of
prime ideals $(p)$ factorizes in, in all the 
cyclotomic subfields of $\mathbb Q(\zeta_n)$.

\begin{Prop}
\label{inequality}
Suppose that $(q,n_1n_2)=1$.\\
{\rm 1)} If $(n_1,n_2)=1$, then $i_q(n_1n_2)\ge i_q(n_1)i_q(n_2)$.\\
{\rm 2)} If $(n_1,n_2)=1$ and $({\rm ord}_q(n_1),{\rm ord}_q(n_2))=1$, then
$i_q(n_1n_2)=i_q(n_1)i_q(n_2)$.\\
{\rm 3)} If $i_q(p_1^{e_1}\cdots p_s^{e_s})$ is prime, then for every $1\le i\le s$ there
exists a $j\ne i$ such that $({\rm ord}_q(p_i^{e_i}),{\rm ord}_q(p_j^{e_j}))>1$.\\
{\rm 4)} We have $i_q(n_1n_2)\ge i_q(n_1)+i_q(n_2)-1$.
\end{Prop}
{\it Proof}. 1+2) If $(d_1,d_2)=1$, then ord$_q(d_1,d_2)={\rm lcm}({\rm ord}_q(d_1),
{\rm ord}_q(d_2))$ and so ord$_q(d_1d_2)\le {\rm ord}_q(d_1) {\rm ord}_q(d_2)$. If
in addition $({\rm ord}_q(n_1),{\rm ord}_q(n_2))=1$, then 
ord$_q(d_1d_2)={\rm ord}_q(d_1) {\rm ord}_q(d_2)$.\\
3) This is a corollary to part 2.\\
4) We have
$$
i_q(n_1n_2)=\sum_{d|n_1n_2}r_q(d)\ge \sum_{d|n_1}r_q(d)+\sum_{d|n_2\atop d\ne 1}r_q(n_1d)
\ge i_q(n_1)+i_q(n_2)-1,
$$
where in the derivation of the latter inequality we used 
part 1 of Lemma \ref{deelbaarheid}.  \qed\\

\noindent The results from \cite{MO} as described above together 
with Lemma \ref{intoirreducibles} (with $q=2$) yield
an explicit formula for $S(n)$.
\begin{Prop} 
\label{Sformula}
Let $i_2$ be as in Lemma {\rm \ref{intoirreducibles}}. We have
$$S(n)=\cases{0 &if {\rm ord}$_2(2n-1)$ is even;\cr
\sqrt{2}^{i_2(2n-1)-1} &if {\rm ord}$_2(2n-1)$ is odd.}$$
\end{Prop}
Note that if ord$_2(2n-1)$ is odd, then for every divisor $d>1$ of $2n-1$ we have
that $\phi(d)/{\rm ord}_2(d)$ is even. Thus $i_2(2n-1)$ is odd and
$S(n)$ is an integer, as a priori it has to be.\\

\noindent The latter proposition together with results from \cite{p1}, then yields the
following first few values of $S(n)$ and $C(n)$.\hfil\break
\vfil\eject

\hfil\break
\centerline{{\bf Table 1:} $S(n)$ and $C(n)$ for small $n$}
\begin{center}
\begin{tabular}{|c|c|c|}\hline

$n$& $S(n)$ & $C(n)$\\ \hline\hline
4&2&1\\ \hline
12&2&1\\ \hline
16&8&2\\ \hline
24&2&1\\ \hline
25&4&1\\ \hline
36&2&1\\ \hline
37&16&2\\ \hline
40&2&1\\ \hline
45&16&2\\ \hline
52&2&1\\ \hline
64&512&30\\ \hline
\end{tabular}
\end{center}

\noindent Using part 4 of Proposition \ref{inequality} and the latter proposition, 
an alternate proof of Corollary \ref{wiseguys} is obtained.
We like to point out that often $S(mn+(m-1)(n-1))>S(m)S(n)$; it can be shown
for example that $S(n^2+(n-1)^2)=S(n)^2>0$ iff ord$_2(2n-1)$ is odd and $2n-1$ is of the form $p^k$
with ord$_2(p)\ne {\rm ord}_2(p^2)$. {}From this equivalence and (\ref{scherper2}), we infer that
$S(n^2+(n-1)^2)\ne S(n)^2$ for almost all $n$ with $S(n)>0$. Part 1 of Proposition \ref{inequality} can also be turned into
a, not so elegant, inequality for $S$.\\

\noindent It is easy to show that ord$_2(p)={\rm ord}_2(p^2)$ iff 
$2^{p-1}\equiv 1({\rm mod~}p^2)$. Primes satisfying the latter congruence
are known as {\it Wieferich primes} and are discussed more extensively in the next section.

\section{On the value distribution of $S(n)$}
By Proposition \ref{Sformula} and the remark that if ord$_2(2n-1)$ is odd, then
$i_2(2n-1)$ is odd, we infer that
$S(n)\in \{0,1,2,4,8,16,32,\ldots \}$.
For $v$ an integer, let $N_v(x)$ denote the number of $n\le x$ for 
which $S(n)=v$ and let $N_v$
be the corresponding set of natural numbers $n$. 
\begin{Prop}
\label{begin}
Let $e\ge 1$ be any natural number. Then $N_{2^e}$ is non-empty.
\end{Prop}
{\it Proof}. Let $p$ be any prime with $p\equiv 3({\rm mod~}4)$,
ord$_2(p)=(p-1)/2$ and $2^{p-1\over 2}\not\equiv 1({\rm mod~}p^2)$ (the prime
$7$ will do), then we claim that
$(p^e+1)/2\in N_e$. Our assumptions imply that $r_2(p)=r_2(p^2)$ and thus,
by Proposition \ref{deelbaarheid}, we infer that $r_2(p^k)=r_2(p)=2$ for $k\ge 1$.
Hence $i_2(p^e)=1+2e$. Note that ord$_2(p^e)$ is odd. By
Proposition \ref{Sformula} it then follows that $S({p^e+1\over 2})=
2^e$. \qed
\begin{Cor}
We have Im$(S)=\{0,1,2,4,8,\ldots\}$.
\end{Cor}
{\it Proof}. This follows from Im$(S)\subseteq \{0,1,2,4,8,\ldots\}$, 
the proposition, $S(1)=1$ and $S(2)=0$. \qed\\

\noindent The sets $N_0$ and $N_1$ are relatively well understood;
we have $N_0(x)=[x]-N(x)$, which together with (\ref{scherper2}) gives
a good estimate and furthermore $N_1=\{1\}$.
It remains to deal with $N_{2^e}$ for $e\ge 1$.
We propose the following strengthening of Proposition \ref{begin}.  
\begin{Conj}
\label{begin+}
Let $e\ge 1$ be any natural number.
There are infinitely many integers $n$ for which $S(n)=2^e$.
\end{Conj}
Put $${\cal P}_m=\Big\{p>2:{p-1\over m}
{\rm ~is~odd~and~ord}_2(p)={p-1\over m}\Big\}$$
and
$${\cal P}'_m=\{p\in {\cal P}_m:
2^{p-1\over m}\not\equiv 1({\rm mod~}p^2)\}.$$
Note that if $p\in {\cal P}_m\backslash {\cal P}'_m$, then ord$_2(p)={\rm ord}_2(p^2)$ and
hence $p$ is a Wieferich prime.
The proof of Proposition \ref{begin} shows that 
the truth of
Conjecture \ref{begin+} follows
if we could prove that
there are infinitely many primes in ${\cal P}'_2$. Although it seems
very likely there are indeed infinitely many such primes, proving this is quite another
matter. A prime being in ${\cal P}_2$ is closely related to
the Artin primitive root conjecture and a prime $p$ satisfying
$2^{p-1\over 2}\equiv 1({\rm mod~}p^2)$, is closely related to
Wieferich's criterion for the first case of Fermat's Last Theorem.
Wieferich proved in 1909 that if there is a non-trivial solution
of $x^p+y^p=z^p$ with $p\nmid xyz$, then $2^{p-1}\equiv 1({\rm mod~}p^2)$ (hence
the terminology Wieferich prime).
Despite extensive computational efforts, the only Wieferich primes ever
found are $1093$ and $3511$ \cite{CDP}. Note that $3511\in {\cal P}_2\backslash {\cal P}'_2$
and that $1093\not\in \cup_{m=1}^{\infty}({\cal P}_m\backslash {\cal P}'_m)$. 
Heuristics suggest that up to $x$ there
are only $O(\log \log x)$ Wieferich primes. If we replace $2$ in
Wieferich's congruence by an arbitrary integer $a$ (i.e., consider primes 
$p$ such that $a^{p-1}\equiv 1({\rm mod~}p^2)$), then it is known
that on average (the $O(\log \log x)$ heuristic
holds true \cite{Mur}, where the averaging is over $a$. Assuming 
the abc-conjecture it is known \cite{Silver}
that there are $\gg \log x/\log \log x$ primes $p\le x$ that are
non-Wieferich primes. For our purposes it would be already enough to 
know that there are only $o(x/\log x)$ Wieferich primes $\le x$, however, even this
is unproved.\\
\indent The situation regarding ${\cal P}_m$ is slightly more promising.
Under GRH we can namely prove
the following result.
Our proof rests on a variation of Hooley's proof 
of Artin's primitive root conjecture that belongs to a class of variations
of this problem dealt with by H.W. Lenstra \cite{L}. This relieves us
from the burden of dealing with its analytic aspects.

\begin{Thm}
\label{stellinkje} {\rm (GRH)}.
Let $m$ be a natural number. Let $\nu_2(m)$ denote the exponent of $2$ in $m$. 
If $m$ is odd or $\nu_2(m)=2$, then ${\cal P}_m$ is empty.
In the remaining cases we have, under GRH,
\begin{equation}
\label{asymptotic} 
{\cal P}_m(x)={2A\epsilon_1(m)\over 3m^2}\prod_{p|m}{p^2-1\over p^2-p-1}{x\over \log x}
+O\left( {x\log \log x \over (\log x)^2}\right),
\end{equation}
where $A=\prod_p\left(1-{1\over p(p-1)}\right)\approx 0.3739558136\ldots$
denotes the Artin constant and
$$\epsilon_1(m)=\cases{1 &if $\nu_{2}(m)=1$;\cr
2 &if $\nu_2(m)\ge 3$.}$$
\end{Thm}
{\it Proof}. If $m$ is odd, then $(p-1)/m$ is even and hence ${\cal P}_m$ is empty. If $4||m$,
then $(p-1)/m$ is odd implies that $p\equiv 5({\rm mod~}8)$ and ord$_2(p-1)=(p-1)/m$ implies
that $({2\over p})=1$. By the supplementary law of quadratic reciprocity the conditions
$({2\over p})=1$ and $p\equiv 5({\rm mod~}8)$ cannot be satisfied at the same time. Hence
${\cal P}_m$ is empty.\\
\indent For the remainder of the proof we assume GRH. Then
it can be shown, cf. \cite{L,W}, that the set of primes $p$ such that ord$_2(p)=(p-1)/m$
satisfies an asymptotic of the form (\ref{asymptotic}) with constant
\begin{equation}
\label{wagstaff}
\sum_{n=1}^{\infty}{\mu(n)\over [\mathbb Q(\zeta_{nm},2^{1/nm}):\mathbb Q]}.
\end{equation}
Let $f$ be an arbitrary natural number. If in 
addition to requiring ord$_2(p)=(p-1)/m$, we also
require $p\equiv 1({\rm mod~}f)$, it is 
readily seen that the sum in (\ref{wagstaff}) has
to be replaced by
$$\delta(f,m)=\sum_{n=1}^{\infty}{\mu(n)\over [\mathbb Q(\zeta_f,\zeta_{nm},2^{1/nm}):\mathbb Q]}.$$
Let $r=\nu_2(m)$. Note that $(p-1)/m$ is odd iff
$p\equiv 1+2^r({\rm mod~}2^{r+1})$. We conclude that (\ref{asymptotic}) holds with
constant $\delta(2^r,m)-\delta(2^{r+1},m)$. It is not difficult to see that
\begin{equation}
\label{tweede}
\delta(2^r,m)-\delta(2^{r+1},m)=\sum_{n=1\atop 2\nmid n}^{\infty}{\mu(n)\over [\mathbb Q(\zeta_{2nm},2^{1/nm}):\mathbb Q]}.
\end{equation}
Let $s|t$. It is
 known (cf. 
\cite[Lemma 1]{M2}) that
\begin{equation}
\label{gradie}
[\mathbb Q(\zeta_t,2^{1/s}):\mathbb Q]=\cases{\varphi(t)s/2 &if $2|s$ and $8|t$;\cr
\varphi(t)s & otherwise.}
\end{equation}
On using the latter 
formula for the degree, we infer after some tedious calculation that the constant in
(\ref{tweede}) equals the constant
in (\ref{asymptotic}). \qed

\begin{Cor}
\label{geenidee}
{\rm (GRH)}. If $e$ is odd or $4|e$, then $N_{2^e}$ is an infinite set.
\end{Cor}
{\it Proof}. If $p\in {\cal P}_{2e}$, then $S({p+1\over 2})=2^{e}$. Now use
that on GRH
the set ${\cal P}_{2e}$ is infinite in case $e$ is odd or $4|e$. \qed\\

\noindent Corollary \ref{geenidee} shows that under GRH Conjecture 1 holds true, provided
that $e$ is odd or $4|e$. It is possible to go further than this, but this requires
a result going beyond Theorem \ref{stellinkje}:
\begin{Thm} {\rm (GRH)}.
\label{stellinkje2} Suppose that $\nu_2(e)\ne 1$.
Let $a$ and $f$ be integers such that
$4e|f$, $\nu_2(4e)=\nu_2(f)$, $(a,f)=1$ and $a\equiv 1+2e({\rm mod~}4e)$. 
There exists an integer $v$ such that for all squarefree $n$ we have
$a\equiv 1({\rm mod~}(f,2en))$ iff $(n,v)=1$. 
The density of primes
$p$ such that $r_2(p)=2e$ and $p\equiv a({\rm mod~}f)$ exists, is
positive, and is given by
$${\epsilon_1(2e)\over \varphi([f,2e])2e}\prod_{p\nmid v}\left(1-{\varphi([f,2e])\over
\varphi([f,2ep])p}\right).$$
\end{Thm}
\begin{Cor} {\rm (GRH)}.
\label{bijnafout}
Suppose that $\nu_2(2e)\ne 2$, $2\nmid f$, $(e,f)=1$ and $(a,f)=1$. Then the
set ${\cal P}_{2e}\cap \{p:p\equiv a({\rm mod~}f)\}$ contains infinitely
many primes.
\end{Cor}
{\it Proof of Theorem} \ref{stellinkje2}. The existence of the density, denoted by
$\delta$, follows
by the work of Lenstra \cite{L}. One obtains that
$$\delta=\sum_{n=1}^{\infty}{\mu(n)c_2(a,f,2en)\over [\mathbb Q(\zeta_f,\zeta_{2en},2^{1/2en}):
\mathbb Q]},$$
with $$c_2(a,f,k)=\cases{1 &if $\sigma_a|_{\mathbb Q(\zeta_f)\cap \mathbb Q(\zeta_{k},2^{1/k})}={\rm id}$;\cr
0 & otherwise},$$
where $\sigma_a$ is the automorphism of $\mathbb Q(\zeta_f):\mathbb Q$ uniquely 
determined by $\sigma_a(\zeta_f)=\zeta_f^a$. Under the conditions of the result one infers that
$$\mathbb Q(\zeta_f)\cap \mathbb Q(\zeta_{2en},2^{1/2en})=\mathbb Q(\zeta_{(f,2en)})$$
and hence $c_2(a,f,2en)=1$ iff $a\equiv 1({\rm mod~}(f,2en))$. Let $n$ be squarefree. Note
that there exists an integer $v$ such that $a\equiv 1({\rm mod~}(f,2en))$ iff $(n,v)=1$.
We infer
that $[\mathbb Q(\zeta_f,\zeta_{2en},2^{1/2en}):\mathbb Q]=
[\mathbb Q(\zeta_{[f,2en]},2^{1/2en}):\mathbb Q]=\varphi([f,2en])2en/\epsilon_1(2e)$, 
by (\ref{gradie}).
We thus find that
$$\delta={\epsilon_1(2e)\over \varphi([f,2e])2e}\sum_{(n,v)=1}{\mu(n)\varphi([f,2e])\over \varphi([f,2en])n}
={\epsilon_1(2e)\over \varphi([f,2e])2e}\prod_{p\nmid v}\left(1-{\varphi([f,2e])\over \varphi([f,2ep])p}\right)>0,$$
where we used that the sum is absolutely convergent and has a summand that is a multiplicative
function in $n$. \qed\\

\noindent {\tt Remark}. An alternative way of proving Theorems \ref{stellinkje} and 
\ref{stellinkje2} is to use the Galois-theoretic method
of Lenstra, Moree and Stevenhagen \cite{LMS}, cf. \cite{S}. This yields, a priori, that, on 
GRH, the
density  is of the form $(1+\prod_{p}E_p)A$, where the $E_p$ are (real) character averages
and hence $-1\le E_p\le 1$ and $E_p=1$ for all but finitely many primes $p$. Moreover, the
$E_p$ are rational numbers and hence the density is a rational multiple of the Artin 
constant $A$.\\

\noindent Now, under GRH, we can prove some further results regarding Conjecture 1.
\begin{Prop} {\rm (GRH)}.
\label{priemniet}
Let $e$ be a natural number. Suppose that $1+2e$ is not a prime number congruent to
$5({\rm mod~}8)$, then $S(n)=2^e$ for infinitely many $n$. 
\end{Prop}
{\it Proof}. If $p\in {\cal P}_{2e}$, then $S({p+1\over 2})=2^{e}$. Hence
the infinitude of ${\cal P}_{2e}$ for $e$ is odd and 
for $4|e$, implies the result in case $1+2e\not \equiv 5({\rm mod~}8)$. 
Next suppose that $1+2e\equiv 5({\rm mod~}8)$ and is not a prime. Then
there exist
natural numbers $e_1$ and $e_2$
with $1+2e=(1+2e_1)(1+2e_2)$ and $1+2e_1\not \equiv 5({\rm mod~}8)$ and $e_2\ge 1$.
Suppose that $n$ is such that
$2n-1=7^{e_2}q$ with 
\begin{equation}
\label{qcondition}
q\not \equiv 1({\rm mod~}6e_1){\rm ~and~}q\not \equiv 1({\rm mod~}14e_1){\rm ~and~}q\in {\cal P}_{2e_1}.
\end{equation} 
The conditions on $q$ ensure that $q\ne 7$ and, moreover,
$(7{\rm ord}_2(7),{\rm ord}_2(q))=1$ and hence $r_2(7^kq)=r_2(7^k)r_2(q)$ for 
every $k$. It follows that
$$i_2(2n-1)=\sum_{d|7^{e_2}}r_2(d)+r_2(q)\sum_{d|7^{e_2}}r_2(d)=(1+2e_1)(1+2e_2)=1+2e,$$
and hence $S(n)=2^e$.
Since $e_1\not \equiv 2({\rm mod~}4)$, it follows by Theorem \ref{stellinkje2} that there are infinitely
many primes $q$ satisfying (\ref{qcondition}). \qed\\

\noindent Some of the cases left open by the previous results are covered by the following result.
\begin{Prop}
\label{priemniet2}
{\rm (GRH)}. Suppose that $1+2e=z$ with $z$ a prime satisfying $z\equiv 5({\rm mod~}104)$, 
then $S(n)=2^e$ for infinitely many $n$.
\end{Prop}
{\it Proof}. The $n$ for which $2n-1=7^2p$, where $({\rm ord}_2(p),21)=3$ and
$p\in {\cal P}_{(z-5)/13}$ will do. Now invoke Theorem \ref{stellinkje2}. \qed\\

\noindent For the values of $S(n)$ hitherto not covered, the following result can
sometimes be applied. Recall that $\omega(n)=\sum_{p|n}1$. We define $\omega_1(n)=\sum_{p||n}1$, i.e. 
the number of distinct prime factors $p$ of $n$ such that $p^2\nmid n$. 
Note that $\omega_1(n)\le \omega(n)$. 
\begin{Thm} 
\label{eenvoorallen}
{\rm (GRH)}. Suppose that $S(n)=2^e$ for some $n$ with
$\omega_1(2n-1)\ge 1$, then
there are infinitely many $n$ for which $S(n)=2^e$. 
\end{Thm}
Our proof rests on the following exchange principle.
\begin{Prop}
\label{i2gelijk}
Let $p$ and $q$ be odd primes and $m$ be a natural number such that 
$(m,2pq)=1$, $r_2(p)=r_2(q)$ and $({p-1\over r_2(p)},{\rm ord}_2(m))=({q-1\over 
r_2(q)},{\rm ord}_2(m))$, then
$i_2(pm)=i_2(qm)$.
\end{Prop}
{\it Proof}. Note that $(p-1)/r_2(p)={\rm ord}_2(p)$. The proof is a consequence of the identity
$$\sum_{d|m}r_2(pd)=\sum_{d|m}r_2(p)r_2(d)({\rm ord}_2(p),{\rm ord}_2(d)),$$
and the observation from elementary number theory that if $a_1,a_2,b$ and $d$ are natural
numbers with $(a_1,b)=d$ and $(a_2,b)=d$, then the 
equality $(a_1,\beta)=(a_2,\beta)$ is satisfied for all $\beta|b$. Hence the conditions
of the proposition imply that 
$$\sum_{d|m}r_2(p)r_2(d)({\rm ord}_2(p),{\rm ord}_2(d))=
\sum_{d|m}r_2(q)r_2(d)({\rm ord}_2(q),{\rm ord}_2(d))=\sum_{d|m}r_2(qd),$$
and thus $i_2(pm)=\sum_{d|m}r_2(d)+\sum_{d|m}r_2(pm)=i_2(qm)$. \qed\\

\noindent {\it Proof of Theorem} \ref{eenvoorallen}. The assumption $\omega_1(2n-1)\ge 1$ implies
that $2n-1=pm$ with $p\nmid m$. By 
Theorem \ref{stellinkje} there exists a number $e_1$ with $\nu_2(2e_1)\ne 2$ such that $p\in {\cal P}_{2e_1}$.
Let $q$ be any prime number such that 
\begin{equation}
\label{nogspecifieker}
q\in {\cal P}_{2e_1}{\rm ~and~}({q-1\over 2e_1},{\rm ord}_2(m))=({p-1\over 2e_1},{\rm ord}_2(m)),
\end{equation}
then by Proposition \ref{i2gelijk}
we infer that $S({qm+1\over 2})=S({pm+1\over 2})=2^e$. By Theorem \ref{stellinkje2} there are infinitely 
many primes $q$ satisfying
(\ref{nogspecifieker}) and thus $S({qm+1\over 2})=2^e$. \qed

\section{Value distribution of $S(n)$ on $n$'s with $2n-1$ squarefree}
Before returning to the 
value distribution of $S(n)$, we address the easier problem of studying the
value distribution of $S(n)$ with $n$ restricted to those
$n$ for which $2n-1$ is squarefree. Let $$V=\{e~:~\mu(2n-1)\ne 0{\rm ~and~}S(n)=2^e{\rm ~for~some~integer~}
n\ge 1\}.$$
(For reasons of space we omit `for some integer $n\ge 1$' in some similar definitions.)
First assume that $e\in V$.
Let
$$r_{\mu}(e)={\rm max}\{\omega(2n-1)~:~\mu(2n-1)\ne 0{\rm~and~}S(n)=2^e \},~{\rm ~and}$$
$${\rm min}_{\mu}(e)={\rm min}\{2n-1~:~\omega(2n-1)=r_{\mu}(e),~\mu(2n-1)\ne 0{\rm~and~}S(n)=2^e \}.$$
In the sequel we will also need the related quantities $r_{\omega_1}(e)$ and
$r_{\omega}(e)$. 
Let $$r_{\omega_1}(e)=\max\{\omega_1(2n-1)~:~S(n)=2^e\}{\rm ~and~}
r_{\omega}(e)=\max\{\omega(2n-1)~:~S(n)=2^e\},$$
and min$_{\omega_1}(e)$ and min$_{\omega}(e)$ be defined as the smallest value of
$2n-1$ for which $\omega_1(2n-1)=r_{\omega_1}(e)$, respectively $\omega(2n-1)=r_{\omega}(e)$ 
and $S(n)=2^e$. In case $e\not\in V$, we put
$r_{\mu}(e)=r_{\omega_1}(e)=r_{\omega}(e)=0$ and leave the associated
minimum quantities undefined.

\begin{Lem}
\label{lemma7}
We have $r_{\mu}(e)\le r_{\omega_1}(e)\le r_{\omega}(e)\le [{\log(2e+1)\over \log 3}]$ for $e\ge 0$.
\end{Lem}
{\it Proof}. If $e\not\in V$, then there is nothing to prove. 
Next suppose $S(n)=2^e$ for some integer $n$. The number $2n-1$ is composed of only primes $p\in {\cal P}$.
Write $2n-1=p_1^{e_1}\cdots p_s^{e_s}$. Then $i_2(2n-1)\ge i_2(p_1\cdots p_s)$. We have
ord$_2(p_i)\le (p_i-1)/2$ and ord$_2(d)\le \phi(d)2^{-\omega(d)}$ for every
divisor $d$ of $2n-1$. Thus
$$i_2(2n-1)\ge i_2(p_1\cdots p_s)\ge \sum_{d|p_1\cdots p_s}2^{\omega(d)}=3^s.$$
If $s>[\log(2e+1)/\log 3]$, then it follows that $i_2(2n-1)>2e+1$ and hence
$S(n)>2^e$. This contradiction shows that 
$s\le [\log(2e+1)/\log 3]$. The proof is concluded on noting that, obviously, 
$r_{\mu}(e)\le r_{\omega_1}(e)\le r_{\omega}(e)$. \qed\\
\indent We now formulate the main result of this section.
\begin{Thm} {\rm (GRH)}.
\label{mainmu}
Given an integer $e$, it is a finite problem to determine whether 
or not it is an element of $V$.
The quantity $r_{\mu}(e)$ can be effectively computed.
\end{Thm}
In order to prove it, it turns out to be fruitful to introduce the
notion of solution tableaux.

\subsection{Solution tableaux}

\noindent Let $s\ge 2$. We define a map $\lambda:\mathbb Z_{>0}^s\rightarrow \mathbb Z_{>0}^s$, 
$(a_1,\ldots,a_s)\rightarrow (l_1,\ldots,l_s)$ as follows. We put
$l_i={\rm lcm}((a_1,a_i),\ldots,(a_{i-1},a_i),(a_{i+1},a_i),\ldots,(a_s,a_i))$ for
$1\le i\le s$. Note that $\nu_p(l_i)=\nu_p(a_i)$ if $\nu_p(a_i)\le \nu_p(a_j)$
for some $j\ne i$ and $\nu_p(l_i)={\rm max}\{\nu_p(a_1),\ldots,\nu_p(a_{i-1}),\nu_p(a_{i+1}),
\ldots,\nu_p(a_s)\}$ otherwise. In particular, if there exist $i$ and $j$ with $i\ne j$
such that
\begin{equation}
\label{gekkie}
\nu_p(a_i)=\nu_p(a_j)={\rm max}\{\nu_p(a_1),\ldots,\nu_p(a_s)\},
\end{equation}
then $\nu_p(l_i)=\nu_p(a_i)$. If (\ref{gekkie}) holds for every prime $p$ (with 
$i$ and $j$ possibly depending on $p$), then 
$(a_1,\ldots,a_s)$ is said to be {\it realizable}. 
Thus in order to determine whether $(a_1,\ldots,a_s)$ is realizable or not, for
each prime divisor $p$ of $a_1\cdots a_s$ one first finds the largest exponent $e_p$
such that $p^{e_p}|a_i$ for some $1\le i\le s$. 
Then one tries to find some $j$ with $j\ne i$ and $1\le j\le s$ such that $p^{e_p}|a_j$.
If this is possible for every $p|a_1\cdots a_s$, then $(a_1,\ldots,a_s)$ is realizable 
and otherwise it is not realizable. The choices of $i$ and $j$ may depend on $p$.
If $(a_1,\ldots,a_s)$ is realizable, 
then $\nu_p(l_i)=\nu_p(a_i)$ for every $1\le i\le s$ and every prime $p$ and hence
$\lambda(a_1,\ldots,a_s)=(a_1,\ldots,a_s)$. The terminology realizable is motivated by
the following result.
\begin{Prop}
We have $(m_1,\ldots,m_s)\in {\rm Im}(\lambda)$ iff $(m_1,\ldots,m_s)$ is realizable.
\end{Prop}
{\it Proof}. `$\Rightarrow$'. By assumption 
$\lambda(a_1,\ldots,a_s)=(m_1,\ldots,m_s)$ for 
some $(a_1,\ldots,a_s)\in \mathbb Z_{>0}^s$. W.l.o.g. assume that
$\nu_p(m_1)={\rm max}\{\nu_p(m_1),\ldots,\nu_p(m_s)\}$. Then it follows that
$p^{\nu_p(m_1)}|(a_1,a_i)|m_i$ for some $1<i\le s$ and so $\nu_p(m_i)=\nu_p(m_1)$.\\
`$\Leftarrow$'. If $(m_1,\ldots,m_s)$ is realizable, then $\lambda(m_1,\ldots,m_s)=
(m_1,\ldots,m_s)$. \qed\\

\noindent Suppose $(m_1,\ldots,m_s)$ is realizable and that $a_i$ is coprime with
$\prod_{j=1}^{i-1}a_j\prod_{j=1}^s m_j$ for $i=1,\ldots,s$, then
$\lambda(a_1m_1,\ldots,a_sm_s)=(m_1,\ldots,m_s)$. Thus if 
$(m_1,\ldots,m_s)$ is in ${\rm Im}(\lambda)$,
then clearly its preimage is an infinite set. However, if
$\lambda(a_1,\ldots,a_s)=(l_1,\ldots,l_s)$, then for $i\ne j$, 
$(a_i,a_j)=(l_i,l_j)$. Thus from a set of merely $s$ numbers all $s(s-1)/2$ gcd's 
$(a_i,a_j)$ can be computed. This is the motivation of introducing the
map $\lambda$ as it leads to the simple expression (\ref{simpleE}) for
$i_2(p_1\cdots p_s)$. Note that alternatively $i_2(p_1\cdots p_s)$ can be
expressed in terms of $r_2(p_i)$, $1\le i\le s$ and the $s(s-1)/2$ gcd's 
$({\rm ord}_2(p_i),{\rm ord}_2(p_j))$.
\begin{Prop}
Given $(a_1,\ldots,a_s)\in \mathbb Z_{>0}^s$, let $\lambda(a_1,\ldots,a_s)=(l_1,\ldots,l_s)$.
We have 
$${\rm lcm}(a_{i_1},\ldots,a_{i_k})={a_{i_1}\ldots a_{i_k}\over l_{i_1}\ldots l_{i_k}}
{\rm lcm}(l_{i_1},\ldots,l_{i_k}),$$
where $k\ge 2$ and $1\le i_1<i_2<\ldots<i_k\le s$.
\end{Prop}
{\it Proof}. It is enough to prove the assertion in case $a_i=p^{e_i}$ with $e_i\ge 0$. W.l.o.g. assume
that $e_1\le e_2\le \ldots \le e_s$. On noting that $l_i=p^{e_i}$ for $i<s$ and that
$l_s=p^{e_{s-1}}$, the result follows after a simple computation. \qed\\

\noindent Let $p_1,\ldots,p_s$ be distinct odd primes. Using the latter proposition we infer
that 
\begin{equation}
\label{simpleE}
i_2(p_1\cdots p_s)=\sum_{v_1=0}^1\ldots \sum_{v_s=0}^1 e_1^{v_1}\cdots e_s^{v_s}
{l_1^{v_1}\cdots l_s^{v_s}\over {\rm lcm}(l_1^{v_1},\ldots,l_s^{v_s})},
\end{equation}
with $(e_1,\ldots,e_s)=(r_2(p_1),\ldots,r_2(p_s))$ and
$\lambda({\rm ord}_2(p_1),\ldots,{\rm ord}_2(p_s))=(l_1,\ldots,l_s)$, where we
order the prime factors $p_1\ldots p_s$ in such a way that $e_i\le e_j$ for $i\le j$
and $l_i\le l_j$ if $e_i=e_j$. We say that $({e_1\atop l_1} \cdots {e_s\atop l_s})$ is the
{\it tableau} associated to $m=p_1\cdots p_s$. We say that 
$({e_1\atop l_1} \cdots {e_s\atop l_s})$ is a {\it solution tableau} if
$\nu_2(e_i)\ne 2$, $e_i$ is even and $l_i$ is odd for every $l_i$ and, moreover, 
$(l_1,\ldots,l_s)$ is realizable. 
The value associated to a solution tableau is the quantity in the right hand
side of (\ref{simpleE}).
The 
following result is a consequence of Theorem \ref{stellinkje2}.
\begin{Prop} 
\label{prop9}
The tableau associated to a squarefree non-prime integer $2n-1$ satisfying
$S(n)=2^e$ is a solution tableau. Given any solution tableau $T$, under GRH, there exist infinitely many 
squarefree integers $2n-1$ such that the associated solution tableau equals $T$.
\end{Prop}
\noindent {\tt Example 1}. Let us consider the tableau $({2\atop 3}{2\atop 15}{8\atop 5})$. Note
that it is a solution tableau corresponding to the value $601$ of $i_2$. Let us try
to find an $m=p_1\cdot p_2\cdot p_3$ having this tableau associated to it. We need to find primes
$p_1,p_2$ in ${\cal P}_2$ and $p_3$ in ${\cal P}_8$ such that 
$({\rm ord}_2(p_1),{\rm ord}_2(p_2))=3$, $({\rm ord}_2(p_1),{\rm ord}_2(p_3))=1$ 
and $({\rm ord}_2(p_2),{\rm ord}_2(p_3))=5$. Any $p_1$ in ${\cal P}_2$ satisfying 
$p_1\equiv 1({\rm mod~}3)$ will fit the bill. E.g. $p_1=79$. Any $p_2$ in ${\cal P}_2$ satisfying
$p_2\equiv 1({\rm mod~}15)$ and $p_2\not\equiv 1({\rm mod~}13)$ will do, 
e.g. $p_2=991$. Finally any $p_3\in {\cal P}_8$ satisfying $p_3\equiv 2({\rm mod~}3)$,
$p_3\equiv 1({\rm mod~}5)$, $p_3\not\equiv 1({\rm mod~}11)$ 
and $p_3\not\equiv 1({\rm mod~}13)$
will do, e.g. $p_3=1721$. We have, as expected, $i_2(79\cdot 991\cdot 1721)=601$. 
(Despite
the numerous conditions on $p_3$, by Corollary \ref{bijnafout} there exist, under GRH, 
infinitely choices for $p_3$.) 
\begin{Prop} 
\label{prop10}
Let $r\ge 1$ be any integer.
The set of solution tableaux associated to the set
$\{m:2\nmid {\rm ord}_2(m),~\mu(m)\ne 0,~i_2(m)\le r\}$ is finite and can
be effectively determined.
\end{Prop}
{\it Proof}. 
Put $\rho=[\log r/\log 3]$. By the proof of Lemma \ref{lemma7} we have 
$\omega(m)\le \rho$. Thus the number of entries in a row
of an associated tableau is bounded by $\rho$. Now fix any $s\le \rho$. Note
that given any integer $k$ there are at most finitely many realizable
$(\lambda_1,\ldots,\lambda_s)$ such that $\lambda_1\cdots \lambda_s/{\rm lcm}(\lambda_1,\ldots,\lambda_s)=k$
(since if $(\lambda_1,\ldots,\lambda_s)$ is realizable, then if $p$ 
divides $\lambda_1\cdots \lambda_s$, it also divides
$\lambda_1\cdots \lambda_s/{\rm lcm}(\lambda_1,\ldots,\lambda_s)$. This observation together with
the remark 
that $i_2(p_1\ldots p_s)\ge e_1\ldots e_s l_1\ldots l_s/{\rm lcm}(l_1,\ldots,l_s)$ shows that there
are only finitely many possibilities for $e_1,\ldots,e_s$, $l_1,\ldots,l_s$ such that
$({e_1\atop l_1}\cdots {e_s\atop l_s})$ is a solution tableau 
and $i_2(p_1\ldots p_s)\le r$. Clearly, these can be effectively determined.\qed\\

\noindent {\tt Example 2}. We try to find all primes $z\equiv 5({\rm mod~}8)$ with $z\le 229$ for
which there are distinct primes $p$ and $q$ both in
$\cal P$ such that $i_2(p\cdot q)=z$. This leads us
to find all $z$ of the above form for which there are $e_1,e_2$ and $w$ satisfying
$1+2e_1+2e_2+4e_1e_2w=z$ with $\nu_2(e_1)\ne 1$, $\nu_2(e_2)\ne 1$ and $w\ge 3$ is odd.
There are solutions precisely for $z\in \{101,157,197,269,317,349,421,509,\ldots\}$. The
associated solution tableau is $({2e_1\atop w}{2e_2\atop w})$. In each case 
one can find $p$ and $q$ corresponding to this solution tableau. E.g., when $z=421$
an associated solution tableau is $({2\atop 5}{38\atop 5})$. A solution corresponding
to this is $p=71$ and $q=174991$.\\

\noindent {\it Proof of Theorem} \ref{mainmu}. 
It is a consequence of Proposition \ref{prop10} that the set of solution tableaux associated
to the set
$$M_e:=\{m:2\nmid {\rm ord}_2(m),~\mu(m)\ne 0,~i_2(m)=2e+1\}$$
can be effectively computed. If this set is empty, then $e\not\in V$. It this set contains
at least one solution tableau, then by Proposition \ref{prop9}, under GRH, it is possible to
find an $m\in M_e$ corresponding to if (cf. Example 1) and thus $e\in V$ by 
Proposition \ref{Sformula}. \qed\\

\noindent We now also have the tools to show that the upper bound in 
Lemma \ref{lemma7} is sharp for infinitely many $e$.
\begin{Prop} {\rm (GRH)}. Let $s\ge 2$. 
If $e=(3^s-1)/2$, then $r_{\mu}(e)=r_{\omega_1}(e)
=r_{\omega}(e)=s$.
\end{Prop}
{\it Proof}. A solution tableau
with $s+1$ columns clearly does not exist. On the other hand 
$T=({2\atop 1}\ldots {2\atop 1})$ ($s$ 
columns) is a solution tableau. Now apply Proposition \ref{prop9}. \qed\\

\noindent {\tt Remark}. If we restrict $e$ to be such that $1+2e$ is prime, then
the upper bound is far from sharp. This is a consequence of part 3 of 
Proposition \ref{inequality}. 

\subsection{The general problem revisited}
Theorem \ref{mainmu} has the following more important variant.
\begin{Thm} {\rm (GRH)}.
\label{main} Assume that the number of primes $p\le x$ such that ord$_2(p)$ is 
odd and $p$ is a Wieferich prime is $o(x/\log x)$.
Given an integer $e$, it is a finite problem to determine whether or not there exist $n$ with
$S(n)=2^e$ and $\omega_1(2n-1)\ge 1$. The 
quantities $r_{\omega_1}(e)$ and $r_{\omega}(e)$ can be effectively computed.
\end{Thm}
For clarity, we first consider some examples. It will be 
convenient to divide the integers $n$ with ord$_2(n)$ odd in
two classes. We say $n=\prod_{i=1}^s{p_i^{e_i}}$ is of type II if there exists $i$ and $j$ such
that $e_i\ge 2$ and $p_i|{\rm ord}_2(p_j)$ and of type I otherwise.\\

\noindent {\tt Example 3}. We try to find all primes $z\equiv 5({\rm mod~}8)$ with $z\le 229$ for
which there are distinct primes $p$ and $q$ both in
$\cal P$ such that $i_2(p^r\cdot q^s)=z$, with $r,s\ge 1$.\\
-Let us assume first that $p^r\cdot q^s$ is of type I. Then we have
to find all $z$ of the above form for which there are $e_1,e_2,r,s$ and $w$ satisfying
$1+2e_1r+2e_2s+4e_1e_2rsw=z$ with $\nu_2(e_1)\ne 1$, $\nu_2(e_2)\ne 1$ and $w\ge 3$ odd.
On making the substitution $f_1=e_1r$ and $f_2=e_2s$, we have to solve
$1+2f_1+2f_2+4f_1f_2w$, with $w\ge 3$ is odd. This equation has the same form as the one arising
in Example 2, except that now there are no restrictions on $f_1$ and $f_2$.
There are solutions precisely for $z\in \{101,{\underbar{109}},157,
{\underbar{173}},197,269,317,349,421,509,\ldots\}$. The two underlined numbers did not arise
in Example 2. A preimage for 109 is readily found using Proposition \ref{priemniet2}. One 
finds that, e.g., $i_2(7^2\cdot 73)=109$. For 173 one finds that the preimage has to be of the form
$p^2\cdot q$ with solution tableau $({2\atop 5}{8\atop 5})$ associated to $p\cdot q$. One finds
that $71^2\cdot 1721$ is the smallest preimage of the required format.\\
-Next assume that $p^r\cdot q^s$ is of type II. In this case $i_2(p^r\cdot q^s)\ge 11+4{\rm min}\{p,q\}$.
We only have to analyze the case where $p\in \{7,23,31,47\}$. The only prime 
$z\equiv 5({\rm mod~}8)$ with $z\le 229$ that is assumed as value turns out to be $197$. This
can happen only if $p$ or $q$ equals 23, e.g., $i_2(23^3\cdot 47)=197$ (we leave the details to
the reader).\\

\noindent {\tt Example 4}. We try to find all primes 
$z\equiv 5({\rm mod~}8)$ with $z\le 229$ for
which there exists an integer $m$ with ord$_2(m)$ is odd and $\omega(m)\ge 3$ 
such that $i_2(m)=z$. In case $\omega(m)\ge 4$ one has
$i_2(m)\ge 2\cdot 2\cdot 2\cdot 2\cdot 3\cdot 5=240$, so assume
that $\omega(m)=3$. In case $m$ is of type I one finds, using
solution tableaux, that there is no such $z$. Indeed, the smallest
prime $z\equiv 5({\rm mod~}8)$ thus occurring is 509. It
arises only if $m=p^2\cdot q\cdot r$ and the solution tableau
associated to $p\cdot q\cdot r$ is $({2\atop 5}{2\atop 5}{2\atop 5})$.
One finds that, e.g., $i_2(71^2\cdot 191\cdot 271)=509$. The smallest
prime $z\equiv 5({\rm mod~}8)$ occurring for a type II integer is
$389$. One has, e.g., $i_2(7^2\cdot 71\cdot 191)=389$.\\

\noindent In Example 3 we managed to reduce the exponents $r$ and $s$ to be $\le 2$ for type I integers. 
The following proposition is also concerned with `exponent reduction'. It is
a slight generalisation of Proposition \ref{i2gelijk} (and so is its proof).
\begin{Prop} 
Suppose that $s\ge 1$. Let $p$ and $q$ be odd primes and $m$ a natural number
such that $(m,2pq)=1$, $(p^{s-1},{\rm ord}_2(m))=1$,
$r_2(q)=\sum_{j=1}^s r_2(p^j)$ (that is $r_2(q)=sr_2(p)$ in 
case $q$ is not a Wieferich prime) and, moreover, 
$({\rm ord}_2(p),{\rm ord}_2(m))=({\rm ord}_2(q),{\rm ord}_2(m))$,
then
$i_2(p^s\cdot m)=i_2(q\cdot m)$.
\end{Prop}
For the prime $q$ above $r_2(q)$ is even, but it is not necessarily the case
that ord$_2(q)$ is odd.\\
\indent  We say that 
$({e_1\atop l_1} \cdots {e_s\atop l_s})$ is a {\it generalized solution tableau} if
$2|e_i$, $2\nmid l_i$ for $1\le i\le s$ and, moreover, 
$(l_1,\ldots,l_s)$ is realizable. In order to find an integer $m$ of type I such
that $i_2(m)=r$ and ord$_2(m)$ is odd, we first determine all generalized solution tableaux
$({f_1\atop l_1}\cdots {f_s\atop l_s})$ associated to $r$. If there are none,
there is no solution. If there is a generalized solution tableau and
$\nu_2(f_i)\ne 2$ for $1\le i\le s$, the generalized solution tableau is 
even a solution
tableau and a squarefree integer $m$ can be found such that $i_2(m)=r$ (under GRH).
W.l.o.g. suppose that $\nu_2(f_i)=2$ for $1\le i\le t$ and $\nu_2(f_i)\ne 2$
for $t<i\le s$. Put $e_i=f_i/2$ for $1\le i\le t$ and $e_i=f_i$ for $t<i\le s$.
Then we want to find an integer $m$ 
of the form $m=p_1^2\cdots p_t^2\cdot p_{t+1}\cdots p_s$, where $p_1,\cdots,p_t$
are non-Wieferich primes and $p_1\cdot p_2\cdots p_s$ has associated solution tableau
$({e_1\atop l_1}\cdots {e_s\atop l_s})$. Under GRH and the assumption that
the number of Wieferich primes with ord$_2(p)$ odd is $o(x/\log x)$, we 
are guaranteed that such an integer $m$ exists. It can be found proceeding as in
Example 1. (Note the `exponent reduction'.)\\
\indent The integers of type II can be dealt with as in Example 3 (and, likewise, there 
are only finitely many possibilities for the smallest prime in the
corresponding number $2n-1$ and each of these cases can be analysed further
using the solution tableau method).\\
\indent This rather informal discussion can be turned into a formal
proof of Theorem \ref{main}. We leave the details to the interested
reader.\\

\noindent The examples together with some additional arguments then yield the
correctness of the following table (for the notation see the beginning of 
Section 6).
\vfil\eject

\hfil\break
\centerline{{\bf Table 2:} Sparse values of $S(n)$}
\begin{center}
\begin{tabular}{|c|c|c|c|c|c|c|}\hline
$1+2e$&${\rm min}_{\mu}(e)$&${\rm min}_{\omega_1}(e)$&${\rm
min}_{\omega}(e)$&$r_{\mu}(e)$&$r_{\omega_1}(e)$&$r_{\omega}(e)$\\
\hline\hline
$5$&$-$&$-$&$7^2$&$0$&0&$1$\\ \hline
$13$&$-$&$-$&$31^2$&0&$0$&$1$\\ \hline
$29$&$-$&$-$&$631^2$&0&$0$&$1$\\ \hline
$37$&$-$&$-$&$127^2$&0&$0$&$1$\\ \hline
$53$&$-$&$-$&$14327^2$&0&$0$&$1$\\ \hline
$61$&$-$&$-$&$3391^2$&0&$0$&$1$\\ \hline
$101$&$7\cdot 631$&$7\cdot 631$&$7\cdot 631$&2&$2$&$2$\\ \hline
$109$&$-$&$7^2\cdot 73$&$7^2\cdot 73$&$0$&1&$2$\\ \hline
$149$&$-$&$-$&$11471^2$&0&$0$&$1$\\ \hline
$157$&$71\cdot 631$&$71\cdot 631$&$71\cdot 631$&2&$2$&$2$\\ \hline
$173$&$-$&$71^2\cdot 1721$&$71^2\cdot 1721$&0&$1$&$2$\\ \hline
$181$&$-$&$-$&$23671^2$&0&$0$&$1$\\ \hline
$197$&$151\cdot 919$&$7^3\cdot 151$&$7^3\cdot 151$&2&$2$&$2$\\ \hline
$229$&$-$&$-$&$248407^2$&0&$0$&$1$\\ \hline
\end{tabular}
\end{center}

\noindent Indeed, using the method of solution tableaux we found that
for $1+2e\le 771$ there exists
an integer $n$ with $\omega_1(n)\ge 1$ and $S(n)=2^e$, with
the exception of the integers $e$ with $1+2e=5,13,29,37,53,61,149,181,229$ and the possible
exception of $461,541$ and $757$. Sometimes the corresponding
values of $n$ were quite large, for example, $i_2(7^{14}\cdot 73)=709$.
Thus merely computing $i_2(n)$ over a large range of $n$ will leave
many small values in the image unreached.  

\section{Analytic aspects}
\noindent We will see that the scarcity of certain values of $S(n)$ is brought
out by analytic number theory.\\ 
\indent Note that $N_2(x)={\cal P}_2(2x-1)$ and hence, by Theorem \ref{stellinkje},
we have under GRH that
$$N_2(x)=A{x\over \log x}+O\left({x\log \log x\over (\log x)^2}\right).$$
Similarly we see that $n\in N_4$ iff $2n-1=p$ with $p\in {\cal P}_4$ or
$2n-1=p^2$ with $p\in {\cal P}'_2$. 
Since ${\cal P}_4$ is empty it follows that $n\in N_4$ iff $2n-1=p^2$ with $p\in {\cal P}_2'$.
On invoking Theorem \ref{stellinkje} we
obtain that, under GRH,
$$N_4(x)\le A{\sqrt{2x}\over \log x}+O\left({\sqrt{x}\log \log x\over (\log x)^2}\right).$$
We actually conjecture that equality holds. Likewise, we conclude that $n\in N_8$ if
$2n-1=p^3$ with $p\in {\cal P}'_2$ or $2n-1=p$ with $p\in {\cal P}_6$. Thus, under GRH,
we deduce on invoking Theorem \ref{stellinkje} that
$$N_8(x)={8Ax\over 45\log x}+O\left({x\log \log x\over (\log x)^2}\right).$$
The asymptotic behaviour of $N_{16}(x)$ is dominated by the number of prime
pairs $(p,q)$ with $p<q$ such that $r_2(p)=r_2(q)=2$, $({\rm ord}_2(p),{\rm ord}_2(q))=1$ 
and $p\cdot q\le 2x-1$.
We are inclined to believe that this number is a positive fraction of all pairs
$(p,q)$ with $p<q$ and such that $p\cdot q\le 2x-1$, which, as is 
well-known \cite{H}, grows asymptotically
as $2x\log \log x/\log x$. Thus we are tempted to conjecture that
$$N_{16}(x)\sim c_0{x\log \log x\over \log x},$$
for some positive constant $c_0$.\\
\indent The following result shows that the values $2^e$ for which 
$r_{\omega_1}(e)=0$ deserve
the predicate {\it sparse}. 
\begin{Prop} 
If $r_{\omega_1}(e)=0$,
then $N_{2^e}(x)\ll \sqrt{x}$.
If $r_{\omega_1}(e)\ge 1$, then
$${x\over \log x}\ll N_{2^e}(x)\le 
(\zeta(2)+o_e(1)){x\over \log x}{(\log \log x)^{r_{\omega_1}(e)-1}\over 
(r_{\omega_1}(e)-1)!},$$ 
where the lower bound holds under the assumption of GRH.
\end{Prop}
{\it Proof}. 
In case $r_{\omega_1}(e)=0$, we are counting a subset of the squarefull 
numbers (a number is squarefull if $p|n$ implies $p^2|n$). As is well-known, 
cf. \cite{Sur}, the number of squarefull integers $m\le x$ grows asymptotically
as $\zeta(3/2)\sqrt{x}/\zeta(3)$. The proof of the lower bound is a 
consequence of Proposition \ref{i2gelijk} and
Theorem \ref{stellinkje2}. As for the upper bound, notice that 
$N_{2^e}(x)$ is bounded above by the number of integers $m\le 2x-1$ such
that $\omega_1(m)\le r_{\omega_1}(e)$. Let $C$ be
an arbitrary positive constant. By the method of 
Sathe-Selberg \cite{H, T} we find that uniformly for $1\le k\le C\log\log x$ we
have 
\begin{equation}
\label{nieuwer}
\#\{m\le x:\omega_1(m)=k\}\sim F({k\over \log \log x})
{x\over \log x}\cdot {(\log \log x)^{k-1}\over (k-1)!},
\end{equation}
where 
$$F(z)={1\over \Gamma(z+1)}\prod_p \left(1+{z\over p}+{1\over p^2-1}\right)
(1-{1\over p})^z,$$
where as usual $\Gamma(z)$ denotes the gamma function. 
{}From this estimate the result follows on noting that $F(0)=\zeta(2)$.\qed\\

\noindent We think the upper bound in the latter result is much closer to the truth:
\begin{Conj} Suppose that $r_{\omega_1}(e)\ge 1$, then
$${x\over \log x}(\log \log x)^{r_{\omega_1}(e)-1}
\ll N_{2^e}(x)\ll {x\over \log x}(\log \log x)^{r_{\omega_1}(e)-1},~x\rightarrow \infty.$$
\end{Conj}

\noindent {\bf Acknowledgement}. We thank Prof. E.W. Clark for pointing out the reference given in
the proof of Lemma \ref{intoirreducibles} and Dr. Breuer and Prof. Shallit for pointing 
out reference \cite{BLM}, respectively \cite{VS}. We thank the referee for attending us on some inaccuracies
in an earlier version.\\
\indent  Part of the
research of the first author was carried out whilst he was temporary assistant professor
in the PIONIER-group of Prof. E. Opdam at the University of Amsterdam.

\medskip\noindent {\footnotesize Max-Planck-Institut f\"ur Mathematik,
Vivatsgasse 7, D-53111 Bonn, Germany.\\
e-mail: {\tt moree@mpim-bonn.mpg.de}}

\medskip\noindent {\footnotesize CNRS-I3S,
ESSI, Route des Colles, 06903 Sophia Antipolis, France.\\
e-mail: {\tt ps@essi.fr}}
\end{document}